\documentstyle[12pt]{article}
\newtheorem{fed}{Definition}[section]
\newtheorem{teo}[fed]{Theorem}
\newtheorem{lem}[fed]{Lemma}
\newtheorem{cor}[fed]{Corollary}
\newtheorem{pro}[fed]{Proposition}
\newtheorem{rem}[fed]{Remark}
\newtheorem{exa}[fed]{Example}
\newtheorem{num}[fed]{}

\def\zC{\hbox{\rm C\hskip -5.4pt\vrule
height 8.0pt width 0.4pt depth -0pt\hskip 4.5pt}}
\def\zR{\ \hbox{\rm R\hskip -11.8pt I \hskip 4.7pt}}

\def\bm{\left(\begin{array}{cc}}
\def\em{\end{array}\right)}
\def\ben{\begin{enumerate}}
\def\een{\end{enumerate}}
\def\hrarr#1#2{\mathrel{ \mathop{\hbox to .5in{\rightarrowfill}} 
\limits^{\scriptstyle#1}_{\scriptstyle#2}  }}
\def\hlarr#1#2{\mathrel{ \mathop{\hbox to .5in{\leftarrowfill}} 
\limits^{\scriptstyle#1}_{\scriptstyle#2}  }}
\def\vdarr#1#2{\llap{$\scriptstyle#1$}\left\downarrow\vcenter to .5in{}\right.\hsize=16cm
\rlap{$\scriptstyle#2$}}\vsize=19cm
\def\ddrarr#1#2{\llap{$\vcenter {\hbox{$\scriptstyle #1$}}$}
\searrow\displaystyle\rlap{$\vcenter {\hbox{$\scriptstyle #2$}}$}}
\def\inc{\subseteq}
\def\noi{\noindent}
\def\csta{C$^*$-algebra}
\def\zC {\hbox{\rm C\hskip -5.4pt\vrule
height 8.0pt width 0.4pt depth -0pt\hskip 4.5pt}}
\def\zR{\ \hbox{\rm R\hskip -8.8pt I \hskip 4.7pt}}

\def\QED{\hskip 0.1truecm\vbox{\hrule\hbox{\vrule\hskip .2truecm
\vbox{\vskip .2truecm}\vrule}
\hrule}\hskip 0.1truecm}
\def\proj{{\cal P} (p)}
\def\projq{{\cal P} (q)}
\def\proo{{\cal P}_{_0} (p)}
\def\prooq{{\cal P}_{_0} (q)}
\def\prof{{\cal P}_{f}(p)}
\def\profq{{\cal P}_{f}(q)}

\def\d+{\Delta^+(p)}
\def\kp{{\cal K}_p }
\def\lp{{\cal L}_p }
\def\ep{{\cal E}_p}
\def\eq{{\cal E}_q}
\def\cu{{\cal U }}
\def\ua{{\cal U} _A}
\def\gla{G_A }
\def\up{\cu (p)}
\def\eps{\varepsilon}
\def\uep{\cu_\eps(A)}

\def\zin{\hbox{\rm $\cap$ \hskip -8.3pt\vrule
height 7.0pt width 0.4pt depth -0pt\hskip 4.5pt}}

\title {\sc \vskip 2truecm Projective spaces of a $C^*$-algebra}
\author
{\normalsize
E. ANDRUCHOW, G. CORACH  AND D. STOJANOFF
\footnote{Partially supported by UBACYT TW49 and TX92, 
PIP 4463 (CONICET) and ANPCYT PICT 97-2259 (Argentina)}{}
}
 
\date{}
\begin{document}

\maketitle
\thispagestyle{empty}
 
\baselineskip=12pt
\begin{quote}
Based on the projective matrix spaces studied by B. Schwarz 
and A. Zaks, we study the notion of projective space associated
to a C$^*$-algebra $A$ with a fixed projection $p$. The resulting space
${\cal P} (p)$ admits a rich geometrical structure as a holomorphic 
manifold and a homogeneous reductive space of the invertible group
of $A$. Moreover, several metrics (chordal, spherical, pseudo-chordal, 
non-Euclidean - in Schwarz-Zaks terminology) are considered, allowing a 
comparison among ${\cal P} (p)$, the Grassmann manifold of $A$ and 
the space of positive elements which are unitary with respect to the
bilinear form induced by the reflection $\varepsilon = 2p-1$. Among 
several metrical results, we prove that geodesics are unique and of minimal
length when measured with the spherical and non-Euclidean metrics.
\end{quote}

\vskip.5truecm
\baselineskip=15pt

\section {Introduction}

There are several papers (\cite{[Ph]}, \cite{[Br]}, \cite{[Z]}) treating the topological and 
metric properties of the space $P= P(A)$ of selfadjoint projections 
of a \csta \ $A$. These properties are used to obtain invariants for the algebra
$A$ (see \cite{[Z]} for instance). Many of these invariants have to do with
problems concerning  the length of curves in $ P$.
There are other papers (\cite{[CPR6]}, \cite{[PR]}, \cite{[W]}) studying $ P$ as a
differentiable manifold (in fact a complemented submanifold of $A$). From
this viewpoint, problems concerning length of curves - e.g.
characterization of curves of minimal length among the curves joining the
same endpoints - can be treated using variational principles, in an
infinite dimensional setting. Since the norms considered in the tangent
bundle of $ P$ do not arise from inner products, this analysis does
not proceed as in the riemannian case, and requires new methods.

In a series of papers \cite{[SZ1]}-\cite{[SZ5]}, B. Schwarz and A. Zaks  
have studied what they call "matrix projective spaces". Their papers inspired
our treatment of the space $P$ as a "one dimensional" projective
space of $A$. Most of the features introduced by Schwarz and 
Zaks for the algebra $M_{2n}(\zC)$ can be carried over in the 
general \csta \ case and we can construct an
identification between the projective space of a \csta  \ and its 
Grassmann manifold of selfadjoint projections.  It should be mentioned that, 
instead of merely
generalizing to the general case the ideas of Schwarz and Zaks, we define a
projective space depending on a fixed projection, which allows us to
deal simultaneously with all "higher dimensional" projective spaces in 
the terminology of Schwarz and Zaks. Many problems
that we study here are not considered in their papers. On the other hand,
many questions treated by Schwarz and Zaks in $M_n(\zC  )$, particularly
those concerning Moebius transformations, will be studied for general
\csta s in a forthcoming paper.

What one gains by taking this standpoint, is the possibility of considering
questions and mathematical objects
related to $ P$, which come up naturally in the projective space
setting, and give interesting information concerning $ P$ and $A$.
Among these, several metrics for $ P$ and the problem of
characterizing their short curves, and the group of projectivities of 
$A$. Moreover, the natural complex structure of the
projective space induces a complex structure on $ P$. Such structure
turns out to be the same that Wilkins obtained by other means in \cite{[W]}.

Let $A$ be a \csta \ and $p \in A$ a projection. Denote by 
$G=G_A$ the group of invertibles of $A$ and 
${\cal U}=\ua $ the unitary group of $A$. 
The orbits $S(p)= \{gpg^{-1} : g \in G_A \}$ and $\cu (p) = \{upu^* : u \in 
\ua \}$ have  rich geometrical and metric properties, studied, 
among others, 
in the papers  \cite{[PR]}, \cite{[W]}, \cite{[CPR1]}, \cite{[CPR2]}, 
\cite{[CPR6]}, \cite{[Br]} and \cite{[Ph]}.
The orbit $\cu (p)$ can be seen as a Grassmann manifold of $A$. 

We denote by $\proj$ the projective space of 
$A$ determined by $p$. It can be defined as the quotient of the set $\kp$ of
partial isometries of $A$ with initial space $p$ by the following 
equivalence relation: two elements $v, w \in \kp$ verify $v \sim_2 w$ if there 
exists $u \in \cu_{pAp}$ such that $v=wu$. In this case we denote by
$[v]=[w]$ the class in $\proj$ (see (\ref{2.9}) for a more 
detailed definition).

The paper \cite{[CPR1]} contains a geometrical study of the set 
$$
S_r=\{ (a,b) \in A\times A : ar=a , rb=b , ba=r \}
$$
where $r$ is an idempotent element of the Banach algebra $A$, and a
corresponding study of the selfadjoint part $R_r$ of $S_r$ if $A$ is a
\csta \ and $r$ is supposed to be selfadjoint. There is an obvious
relation between the spaces $\kp$ and $S_r$, and this paper may be
seen as a kind of continuation of \cite{[CPR1]}. Many constructions done in this
paper can be generalized to the Banach algebra setting. We choose the
\csta \ case in order to keep the paper into a reasonable size.

We define a natural C$^\infty$ manifold structure on $\proj$, 
and the chordal and spherical metrics
generalizing \cite{[SZ2]}.  We show (Theorem \ref{3.5}) that the 
spherical metric has curves of minimal length, 
which are in fact the geodesics determined by the C$^\infty$ 
homogeneous reductive structure induced on $\proj$ 
by the natural action of $\ua$,  given by left multiplication. 

We show that the projective space $\proj$ is diffeomorphic to the Grassmann
manifold $\ep = \{$ projections $ q \in A : q \sim p \}$, where $\sim$ denotes
the usual equivalence of projections (see Theorem \ref{2.14}). Via this 
diffeomorphism, we charac\-terize the chordal metric
of $\proj$ as the metric induced by the norm on $\ep$. Also the 
spherical metric of $\proj$ is identified with the geodesic metric defined in 
$\ep$ by its natural Finsler structure (see (\ref{2.17})). 
Note that $\ep$ is a discrete union
of unitary orbits of projections of $A$. Then $\ep$ has the same local 
geometrical structure as $\cu (p)$. We show (Proposition \ref{3.7}) that there 
exists a unique geodesic
of minimal length joining any two points of $\proj$ which have spherical 
distance less than $\pi /2$. This result was unknown for the 
Grassmann manifolds.

We define the group of projectivities of $\proj$, using the action of $\gla$ 
on $\proj$ by left multiplication. The set of ``finite points'' $\prof$ 
of $\proj$ is characterized in terms of the chordal 
and spherical metrics. For example, 
it is shown that $\prof$ is exactly the set of points whose spherical
distance to $[p]$ is less than $\pi /2$ (see Proposition \ref{4.7}). 
A consequence of this fact is that $\prof$ is homeomorphic to the linear 
manifold $H_p = (1-p)Ap$ (see Proposition \ref{4.7}). The Moebius 
maps are defined and their domains are characterized (\ref{4.8}). 
They are of particular interest in the  case of the algebra
$A=B^{2\times2}$ for a \csta \ $B$, when $p = \bm  1&0\\0&0
\em $ (see (\ref{4.5})).

A holomorphic structure is defined on $\proj$ (Theorem \ref{5.7}) 
via the local homeorphisms
mentioned before. Also a homogeneous reductive structure is introduced, using 
the natural action of the Lie group $\gla$ given by the projetivities 
(see (\ref{5.7})). 

Finally we consider the pseudo-chordal and non-Euclidean metrics (generalizing
the definitions of \cite{[SZ2]}) on the unit disc $\d+$, defined as the orbit of 
$[p]$ in $\proj$ by the action of the group of $\eps$-unitaries $\uep$ by 
left multiplication,  where $\eps$ is the symmetry $\eps = 2p-1$.
The disc $\d+$ is characterized in several ways 
(Propositions \ref{6.9} and \ref{6.12}) 
and the pseudo-chordal and non-Euclidean metrics are showed to be the 
translation of the natural metrics of the space $\uep ^+ = A^+ \cap \uep$ studied 
in \cite{[CPR4]}, \cite{[CPR5]} and \cite{[CPR6]}  (see Theorem \ref{6.17}). 
Also a C$^\infty$ manifold structure is defined on $\d+$ and the natural action of 
$\uep$ converts $\d+$ a homogeneous reductive space with a Finsler metric. 
We show that the geodesics become curves of minimal length, and therefore the 
non-Euclidean metric is rectifiable.

The connections between the work by L. G. Brown, G. K. Pederesen \cite{BP1}, 
\cite{BP2}, G. K. Pedersen \cite{Ped} and S. Zhang \cite{[Z]} and our work
is not completely understood, but they may have interesting consequences
in our context. We expect to deal with these matters in a future work. However,
we should warm the reader about the completely different use of the word 
"projective" done by us and by those who deal with projective \csta s.

\section{ The projective space.}

Let $A$ be a \csta , $\gla$ the group of invertibles of $A$ and 
$\ua$ the unitary group of $A$. 
Let $p=p^2=p^* \in A$ be a fixed projection. If $C$ is a
subset of $A$, $Cp$ denotes the set $\{cp : c \in C\}$.

Usually, one regards the space of projections $\ep = \{q : q\sim p\}$
as an homoge\-neous space (i.e. quotient of) the unitary group of $A$. Here we propose
an alternate view of $\ep $, considering another natural action,
of the analytic Lie group $G_A$. Using the fixed projection $p$, one can regard
the elements of $A$ as $2\times 2$ matrices. We shall consider the set of matrices
with second column equal to zero, and introduce there a equivalence relation.
It will be readily clear that $G_A$ acts on the quotient $\proj$
(by left multiplication), and that this space $\proj$ is 
homeomorphic to $\ep $.

\begin{fed}\label{2.1} \rm
Let $A$ be a \csta \ and $p\in A$ 
a projection. We consider the following subsets of $A$
$$
\lp  =\{a \in Ap : \hbox{ there exists } b \in pA \hbox{  with }
ba = p \}
$$
and 
$$
\kp =\{v  \in Ap : v^*v=p \}.$$
Note that $\kp $ consists of the partial isometries of $A$ with
initial space $p$.
\end{fed}

\begin{rem} \label{2.2} 
\rm If $A=M_n (\zC )$ and $p\in A$ is a 
projection, then 
$$
\lp \ =\gla \cdot  
p\ =\ \{ \ a \in Ap :  \ rank(a)=rank(p) \ \}.
$$
Analogously $\kp =\ua \cdot  p$. In general, the inclusions 
$\gla \cdot  p \subset \lp $ and $\ua \cdot  p \subset \kp $ are strict
(e.g. $A= {\cal B}(H)$).
\end{rem}

\begin{fed}\label{2.3} \rm Let $p\in A$ 
a projection. Let us state the following equivalence relations:
\ben
\item The relation $\sim _1$ in $\lp$ : $a_1p \sim _1 a_2p$ 
if there exists $h \in G_{pAp}$ such that $a_1ph=a_2p$.

\item The relation $\sim _2$ in $\kp$ : $v_1p \sim_2 v_2p$ 
if there exists $w \in {\cal U}_{pAp}$ such that $v_1pw=v_2p$.
\een
\end{fed}

\begin{rem}\label{2.4}
\rm  If we write the elements of $A$ as $2\times 2$ matrices 
using $p$, then
$$
\lp =\{ \bm a & 0 \\ b & 0 \em \in A :
\hbox{ there exists } 
\bm c & d \\ 0 & 0 \em \in A 
\hbox{ such that }  \bm ca +db & 0 \\ 0 &
0 \em = p \}
$$
with an analogous description for $\kp $. The equivalence relation
$\sim_1$ is given by
\begin{equation}\label{2.5}
\bm a_1&0 \\ b_1&0 \em \sim_1 \bm a_2 & 0
\\ b_2 & 0 \em \ \hbox{ if } \ 
\bm a_1&0 \\ b_1&0 \em = \bm a_2h & 0
\\ b_2h & 0 \em , \ 
 \end{equation}
for some $ h \in G_{pAp}$ and analogously for $\sim_2$. 
Note that $\sim_2$ is just the restriction
of $\sim_1$ to $\kp $.
\end{rem}

\begin{rem}\label{2.6}
\rm It is easy to see that, for $a \in Ap$,
$$
a \in \lp \Leftrightarrow a^*a \in G_{pAp}.
$$
Therefore, if $a \in \lp$, the unitary part of the right polar
decomposition of $a$ is  $u = a |a|^{-1} \in A$, where $|a|^{-1}$ is the
inverse of $|a| = (a^*a)^{1/2}$ in $pAp \inc A$. Note that, by
construction, one gets that $a \sim_1 u$ and $u \in \kp$.
\end{rem}

\begin{cor}\label{2.7} If $a \in \lp $, then there
exists $u\in \kp $ such that $a \sim_1 u $ in $\lp $.
\end{cor}

\begin{cor}\label{2.8} If $g \in Gp$, then there exists $v
\in {\cal U}$  such that $gp \sim_1 vp$.
\end{cor}
 
\noi Proof. \rm Suppose that $A$ is faithfully represented in a Hilbert
space $H$. Then  $gp\in \lp $ and $g(1-p)\in {\cal L} _{1-p}(A)$ and 
therefore there exist partial isometries $v_1,v_2 \in A$
$$
v_1 :p \to Im(gp)=M \, \hbox{ and } \, v_2 :1-p \to Im(g(1-p))=N ,
$$
such that $gp=v_1 |gp|$, $g(1-p)= v_2 |g(1-p)|$, $|qp| \in G_{pAp}$ and
$|g(1-p)| \in G_{(1-p)A(1-p)}$. Since $g \in G$, then $M \oplus
N=H$. Let $q_1$ be the
orthogonal projection onto $M$ and $q_2$ be the orthogonal projection onto
$N$. Then $q_1 = v_1v_1^* \in A$ and $q_2 = v_2 v_2^* \in A$. Moreover, 
it is easy to see that  that $\|q_1 q_2\| <1$. 
Hence $\|q_2q_1q_2\|<1$, and $q_2 -q_2q_1q_2=q_2(1-q_1)q_2 \in G_{q_2Aq_2}$. 
Therefore $(1-q_1)q_2\in {\cal L} _{q_2}(A)$ and its
polar decomposition is $(1-q_1)q_2= u|(1-q_1)q_2|$, 
with $u \in A$  a partial isometry
$u:q_2 \to 1-q_1$. 

Let $v_3=uv_2 \in A$. Then $v_3$ is a partial isometry from $1-p$ to $1-q_1$, 
$v=v_1+v_3$ is a unitary element of $A$ and 
$gp=v_1|gp|=vp|gp|$, hence $vp \sim_1 gp$. \QED

\begin{fed}\label{2.9}  \rm Let $A$ be a \csta \ and $p\in A$ 
a projection. We define  the projective space of $A$ determined by $p$:
$$
\proj  =\lp /\sim_1 \quad \hbox{ and } \quad 
\proo =Gp/ \sim_1.
$$ 
\end{fed}

\begin{rem}\label{2.10} \rm  The previous results prove 
that the inclusion map  
$\kp  \hookrightarrow \lp $ induces the bijection
$$
\kp / \sim_2 \ \to \ \lp /\sim_1 = \proj .
$$ 
Analogously the inclusion map  $\ua \cdot  p   
\hookrightarrow \gla \cdot  p $
induces the bijection
$$
\ua \cdot  p / \sim_2 \ \to \ \gla \cdot  p /\sim_1 = \proo .
$$
In both sets we shall consider the quotient topology induced by the norm
topology of $A$. It will be shown that these bijections are in fact homeomorphisms.
\end{rem}

\begin{fed}\label{2.11} \rm 
Denote $\ep  =\{ q \in A : q^* = q =q^2 \hbox{ and } q \sim p \}$, 
where $\sim$
denotes the (Murray-von Neumann) equivalence relation for 
projections of a \csta ,
i.e. for two projections $p,q \in A$,  $p \sim q$ 
means that there exist $v \in A$ such that $vv^*=q$ and $v^*v=p$.
\end{fed}

\begin{rem}\label{2.12}\rm \ 

\ben
\item The group $G_A$ acts on $\proj $ and 
$\proo $ by left multiplication. Namely, if $g \in G_A$ and $[a] \in \proj $,
put $g \times [a] = [ga]$. The same definition
works in $\proo $. Occasionally, we shall consider the restriction of
this action to $\ua$.

\item $\ua $ acts also on $\ep $, by means of $(u, q) \mapsto uqu^*$.
The orbits of this action lie at distance 
greater or equal than $1$ (computed with the norm of $A$) - it is a standard fact that projections at
distance less than $1$ are unitarily equivalent with a unitary element
in the connected component of 1 (see \cite{[RS-N]}, for example). Therefore each one of
these orbits consists of a discrete union of connected components of the space of
projections of $A$ (also called the Grassmannians of $A$). 
These  are well known spaces, which have rich
geometric structure as homogeneous reductive spaces and C$^\infty$
submanifolds of $A$ (see \cite{[PR]} and \cite{[CPR6]}).

Therefore $\ep $ is a submanifold of the Grassmannians of $A$. If,
additionally, $\ua $ is connected, then each component of $\ep $ is the
unitary orbit of a projection.
\een
\end{rem}

\noi We shall see that the space $\proj $ endowed with the quotient
topology is homeomorphic to $\ep $, therefore inheriting the
differentiable structure of the Grassmannians. 
The mapping $\kp  \to \ep $ given by $ v \mapsto vv^*$ is clearly continuous 
and surjective. Clearly it defines a continuous surjective map 
\begin{equation}\label{2.13}
\varrho_p: \proj  \to \ep  \quad \varrho_p([v])=vv^*.
\end{equation}

\begin{teo}\label{2.14} Let $A$ be a \csta \ and $p\in A$ 
a projection. Then $\varrho_p$  is a homeomorphism.
Moreover, if $[v] \in  \proj $ and $q=\varrho_p ([v])$, then the following
diagram commutes:
$${
\begin{array}{ccc}
\ua  & \hrarr{\pi_{[v]}}{} & \proj    \\
           & \ddrarr {\pi_q }{} &  \vdarr {}{\varrho_p} \\
           &                   & \ep    
\end{array} } 
$$
where $\pi_{[v]} (u)= [uv]$ and $\pi_q (u)= uqu^*$, for $u \in \ua $.
\end{teo}

\noi Proof. \rm Let us to prove that $\varrho_p$ is one
to one. Suppose that  $v_1,v_2 \in \kp $ with $v_1 v_1^*=v_2 v_2^*$ and  let 
$w=v_2^*v_1$. Note that $w \in {\cal U}_{pAp}$ and that 
$$
v_2 w=v_2 v_2^*v_1= v_1 v_1^* v_1= v_1 p=v_1 ,
$$
that is, $v_1 \sim_2 v_2$. Straightforwrad computations show that the diagram commute.
The map $\pi_{[v]}$ is continuous and $\pi_q$ has continuous local
cross sections (see \cite{[PR]}). 
Using the diagram, these facts imply that $\varrho_p $ is an
open mapping.\QED
\medskip

At the beginning of the section we noted that the sets $\lp /\sim_1$
and $\kp / \sim_2$ are in a bijective correspondence. Now we shall prove that their
respective quotient topologies also coincide.

\begin{pro}\label{2.15}If $\lp /\sim_1$
and $\kp / \sim_2$ are endowed with their quotient topologies, then
the inclusion map  $\kp  \hookrightarrow \lp $
induces the homeomorphism
$$
\kp / \sim_2 \ \to \ \lp /\sim_1 .
$$
\end{pro}

\noi Proof. \rm Suppose that the algebra $A$ is represented on a 
Hilbart space $H$. It suffices to prove that the mapping
$$
\lp / \sim_1 \ \to \ \ep  \, \hbox{ given by } \, [a] \mapsto P_{a(H)} \quad , 
\quad a \in \lp ,
$$
is continuous, where $P_{a(H)}$ denotes the projection onto the range of $a$.
As shown before, $a \in \lp  $ implies $|a| \in G_{pAp}$. Now, if $(a^*a)^{-1}$
is the inverse of $a^*a$ in $pAp$, then 	
$P_{a(H)} =a (a^*a)^{-1}a^*$, since $a|a|^{-1}$ is a partial isometry with
initial space $p$ and final space $P_{a(H)}$. The result follows
reasoning as in the previous theorem.\QED

\begin{cor}\label{2.16} Let $A$ be a \csta , $p \in A$ 
a projection and  $[a] \in \proj$.
\ben
\item The  orbits of $[a]$ 
by the action of $\ua$ and $\gla$ coincide. That is, 
$$
\cu_{[a]} := \{[ua]: u \in \ua \}=\{[ga] : g \in \gla \} := S_{[a]}.
$$
\item  The connected component of $[a]$ 
in  $\proj$ is contained in $S_{[a]}$.
\item If $G_A$ (or equivalently, $\ua$) is connected, then the connected
component of $[a]$
in  $\proj$ is exactly $S_{[a]}$.
\een
\end{cor}
 
\noi Proof. 
It is well known (see \cite{[PR]} or \cite{[CPR6]}) that 2 and 3 
are true in $\ep$. So they are also true in $\proj $ using the homeomorphism
$\varrho _p$. We know that 1 is true for $a \sim_1 p$, by (\ref{2.10}). 
For any other $[a] \in \proj$  denote by $q = \varrho_p([a])$. Then 
the result follows applying (\ref{2.10}) to $\projq $.\QED

\begin{rem}\label{2.17}\rm 
$\proj $ has C$^\infty$ 
differentiable, homogeneous and (unitary) reductive 
structure induced by the homeomorphism with the space
$\ep $ which, as pointed out before, has rich geometric structure
studied in  \cite{[PR]}  and  \cite{ [CPR6]}. Let us recall some facts:

\ben
\item  The space of projections of $A$, with the norm
topology is a discrete union of unitary orbits of projections. Each orbit
is a C$^\infty$ submanifold of $A$, and a C$^\infty$ homogeneous space
under the action of Lie-Banach group $\ua$. 
The tangent space at a
given projection $p$ identifies with the $2 \times 2$ matrices (in terms
of $p$) which are selfadjoint and have zeros in the diagonal.

\item The space of projections of $A$
admits a natural reductive structure, which induces a linear
connection. The invariants of this connection can be explicitely computed.
It is torsion free and the curvature tensor is given by
$$ 
R(x,y)z=\Big{[} \ [z, p\, ] ,\ [x,y] \ \Big{]} ,
$$
where $[a,b]=ab-ba$ for $a,b \in A$.

\item The geodesic curves of this connection can be computed.
The unique geodesic $\gamma$ with $\gamma (0)=p$ and $\dot
\gamma (0)=x$ is given by
$$
\gamma (t) = e^{t[x,p\, ]}p\, e^{-t[x,p\, ]}.
$$

\item  There is a natural invariant Finsler metric on the space of
projections of $A$, namely the usual norm of $A$ in every tangent space. This
metric has remarkable properties which will be recalled
later.
\een
\end{rem}

\section{The chordal and spherical metrics on $\proj $.}

In this section we introduce the two metrics on
$\proj$ referred in the title. They are the operator theoretic analogues of
the metrics considered in \cite{[SZ2]} for projective (finite dimensional) matrix
spaces.

\begin{fed}\label{3.1} \rm 
If $[a],[b] \in \proj $  for $a,b \in \kp $, the
chordal distance between $[a]$ and $[b]$ is
$$
d_c ([a],[b])= \|\varrho_p ([a])- \varrho_p ([b])\|= \|aa^* - bb^*\|.
$$
\end{fed}

\begin{rem}\label{3.2}
\rm  

\ben
\item This is the metric given by the natural (norm) metric of $\ep $, 
and therefore induces the already considered quotient topology on
$\proj $.
\item If two elements lie in different connected components of $\proj $, 
then their chordal distance is greater or equal than 1.
\item This metric is invariant under the action of $\ua$.
\item If $a\in \kp $, $u \in \ua$  and $q=aa^*$, then
$$
d_c ([ua],[a])= \|uq -qu\|= \max \{\|qu(1-q)\|, \|(1-q)uq\|\} \le 1.
$$
In particular, this shows that our definition agrees with the chordal
distance considered in \cite{[SZ2]} for the case $A=M_n(\zC )$.
\een
\end{rem}
\medskip
By means of the chordal metric we can compute length of curves in $\proj$. 
Given a curve $\gamma : [0,1] \to \proj$ consider the length 
$\ell (\gamma )$ as 
$$\ell (\gamma ) =
\sup _\pi \sum _{i=1}^n \ d_c ( \gamma (t_i) , \gamma (t_{i+1} ) ), 
$$
where the supremum is taken over all partitions $\pi$ of $[0,1]$.

\begin{fed}\label{3.3} \rm 
If $[a], [b] \in \proj $ lie in
the same connected component, then
$$
d_r ([a],[b])= \inf \{ \ell (\gamma): \ \ 
\gamma : [0,1] \to \proj \hbox{ with } \gamma(0)=[a] 
\hbox{ and } \gamma(1)=[b]\}, 
$$
where the infimum
is taken over all rectifiable curves (i.e. $\ell (\gamma) < \infty$). 
Define $d_r ([a],[b])=
\infty$ if $[a]$ and $[b]$ lie in different connected components. $d_r$ is 
called the rectifiable metric.
\end{fed}

\begin{rem}\label{3.4}\rm
The differentiable structure of $\proj $ allows one to compute $d_r$ using
 C$^1$ curves. In this case 
$$
\ell (\gamma)= \int_0^1 \|\dot \gamma (t) \| dt.
$$
This fact means that $d_r$ is the translation of the 
rectifiable (or geodesic) metric in $\ep $
by means of the diffeomorphism $\varrho _p$. This metric has been extensively
studied (\cite{[PR]}, \cite{[CPR6]}, \cite{[Br]}, \cite{[Ph]}).
The following theorem collects some of its properties:
\end{rem}

\begin{teo}\label{3.5} Let $p \in A$ a projection and $q \in
\cu (p)$. Then 
\ben
\item $\|p-q\| <1$ if and only if $d_r (p,q) < \pi/2$.

\noi In this case: 
\item there exists a geodesic in $ \cu (p) $
joining $p$ with $q$ whose length is minimal and therefore equals $d_r (p,q)$.
\item $d_r (p,q)=$  arcsin$(\|p-q\|)$. 
\een
\end{teo}

\noi Proof.  We use a result of  \cite{[PR]} which says that if $\|p-q\|<1$ then there
exists a geodesic curve $\gamma$ joining $p$ and $q$ with $\ell
(\gamma)=d_r (p,q)\le \pi/2$. 
Let $x$ be the velocity vector of this geodesic. That
is, $\gamma (t)= e^{tx}pe^{-tx}$.  In matrix form (in terms of $p$):
$$
x=\bm 0&-a^* \\ a &0 \em 
$$
and $\ell (\gamma)=\|x\|=\|a\|\le \pi/2$.

On the other hand  $\|p-q\|=\|e^xp-pe^x\|$. Easy calculations show that 
$$
e^x =\bm \cos(|a|) & -a^* f(|a^*|) \\ a f(|a|) & \cos(|a^*|)
\em 
$$
where $f(t)= \frac{\sin (t)}{t}$, defined for $t\ge0$. It is easy to see
that $\|a f(|a|)\| = \| \sin (|a|) \|$ and the same for $a^*$.
Since  $\|a\|\le \pi /2 $, we can deduce that $ \| \sin (|a|) \|=
 \sin (\|a\|)$. Therefore 
$$
\|p-q\|=\|pe^x-e^xp\|=\max \{\|a f(|a|)\|,\|a^* f(|a^*|)\|\}=\sin (\|a\|).
$$
This shows one implication in 1 via the formula 3.  
Now, if $d_r(p,q)<\pi /2$, the argument consists on taking a short curve
$\rho$ joining $p$ and $q$ and a partition such that each pair of contiguous
projections have chordal distance less than one. By the previous result 
a polygonal $\gamma$ of geodesics shorter 
than $\rho$ can be constructed. Therefore
the sum of its lengths (= $\ell (\gamma )$) 
is less than $\pi/2$. Finally we use the following result (Lemma 3 of \cite{[Br]}):
Given three projections $r,s,w$ and nonnegative numbers $t_1, t_2$ such
that $t_1 +t_2 < \pi/2 $, $\|r-s\|= \sin t_1 $ and $\|s-w\|= \sin t_2$, then
$\|r-w\|\le \sin (t_1 + t_2)$.
In our case, it can be easily deduced that 
$\|p-q\| \le \sin (\ell (\gamma )) < 1 $. This concludes the proof.\QED

\begin{rem}\label{3.6}\rm 
The contents of the theorem above are essentially known.
Items 2 and a part of 3 of (\ref{3.5}) were proved in \cite{[PR]}. 
Phillips \cite{[Ph]} proved  equality 3 and the existence of curves of
minimal length.
The proof of 3 presented here uses the original ideas of 
Porta and Recht \cite{[PR]}. It should be noted that in \cite{[PR]} 
the results are stated
in terms of selfadjoint symmetries (i.e. elements $\varepsilon \in A$ with 
$\varepsilon^2=1$ and
$\varepsilon^*=\varepsilon$). One passes from projections to symmetries by 
$p \to \varepsilon=2p-1$, and therefore the metric in the space of 
symmetries carries a factor $2$ with respect to the metric in the space of 
projections.
\end{rem}

\medskip

In \cite{[PR]} it was shown that if $\|p-q\|<1/2$ then $p$ and $q$ are joined
by a {\bf unique} geodesic of the linear connection. The following 
proposition says that this fact still holds if $\|p-q\|<1$. Note 
that if $\|p-q\|<1$ there are many curves joining $p$ and $q$
with minimal length. What the following statement says is that
only one of them is a geodesic.

\begin{pro}\label{3.7} Let $D=\{ z \in A : z^*=-z 
\, , \, pz=z(1-p) \hbox{ and } \|z\|<\pi/2 \}$. Then 
$$
\exp : D \to \{q \in \ep  : \|p-q\| <1\} , \quad \exp (z)=e^z p e^{-z}
$$
is a C$^\infty$ diffeomorphism.
\end{pro}

\noi Proof.  Let $x \in A$ with $x^*=-x$ and 
$px=x(1-p)$. As in the previous result
$$
x=\bm 0&-a^* \\ a &0 \em 
$$
 and
$$
e^x =\bm \cos(|a|) & -a^* f(|a^*|) \\ a f(|a|) & \cos(|a^*|)
\em 
$$
with  $f(t)$ as in 3.5. Put $\varepsilon=2p-1$. Then condition 
$pz=z(1-p)$ becomes $\varepsilon 
z=-z\varepsilon$, and therefore $e^z\varepsilon =\varepsilon e^{-z}$.

Clearly 3.5 implies that the mapping $\exp$ is surjective between the
domains considered. Let us prove that it is also one to one. Suppose that 
$z_1,z_2 \in D$ satisfy  $e^{z_1}\varepsilon e^{-z_1}=e^{z_2}\varepsilon
e^{-z_2}$. Then $e^{2z_1}\varepsilon 
= e^{2z_2}\varepsilon$, and since $\varepsilon$ is invertible, this 
implies $e^{2z_1}=e^{2z_2}$. Both exponentials have matrix forms as above, 
$$
e^{2z_i} =\bm \cos(|a_i|) & -a_i^* f(|a_i^*|) \\ a_i
f(|a_i|) & \cos(|a_i^*|)
\em 
$$
with $\|a_i\|<\pi$, $i=1,2$. The function $\cos$ is a diffeomorphism on the 
set of positive elements of $A$ with norm strictly less than $\pi$. Therefore 
$\cos (|a_1|)=\cos(|a_2|)$ implies $|a_1|=|a_2|$. Another functional calculus 
argument shows that both $f(|a_i|)$, $i=1,2$ are invertible elements of $A$, 
and therefore $a_1=a_2$. Moreover, this same sort of argument shows that 
$\exp$ is a diffeomorphism, since its inverse can be explicitely computed.
\QED

\begin{rem}\label{3.8}\rm 
Equality 3 of (\ref{3.5}) implies that the metric
$d_r$ on $\proj $ agrees with the spherical distance defined in \cite{[SZ2]} as the
arcsin of the chordal distance, under the hypothesis that the
chordal distance between the pair of points is less than one.
Moreover, if the diameter of $\cu (p)$ (using the metric $d_r$) is $\pi/2$, 
then the spherical distance equals $d_r$ in each unitary orbit of $\proj$. 
Indeed, easy computations show that if $q,r \in \cu (p)$ and 
$d_r(q,r) = \pi/2$ then $\|q-r\|=1$. In  \cite{[Ph]} it is shown that  a large 
class of  \csta s satisfy this diameter condition.
\end{rem}

\section{Projectivities and finite points of $\proj$.}

Let $H$ be a Hilbert space,  $A \inc {\cal B} (H)$ a \csta \ and
$p \in A$ a projection.

\begin{fed}\label{4.1} \rm
Let $g \in \gla$. We denote by $T_g :\proj \to
\proj $  the map
$$
T_g ([a]) = [ga] \quad , \quad [a] \in \proj .
$$
It is clearly well defined and is a diffeomorphism. 
Following \cite{[SZ2]}, we call these maps the {\bf projectivities} of $\proj$.
\end{fed}

\begin{num}\label{4.2} \rm If 
we identify $[u] \in \proj $ for $u \in \kp $
with $q = uu^* \in \ep $, we can describe the projection in
$\ep$ which corresponds to $T_g ([u])$.
Note that if $a \in \lp $ and $u \in \kp $
is the partial isometry
appearing in the polar decomposition of $a$ in $A$,
then $[a] = [u] \in \proj $. Since $a(H) = u(H)$ by
construction, we deduce that $q = uu^* = P_{a(H)}$.
The same happens for $ga \in \lp $. Therefore
$T_g([a]) = [ga]$ can be identified with the
projection $T_g(q) = P_{ga(H)}$. Note also that $gqg^{-1}(H) =
ga(H)$. Therefore $T_g(q) $ is the projection onto the
image of the idempotent $gqg^{-1}$. Therefore (see \cite{[CPR6]}),
$$
\begin{array}{rl}
T_g(q) = P_{gqg^{-1}(H)}& =
gqg^{-1}(gqg^{-1})^*(1+(gqg^{-1}-(gqg^{-1})^*)^2)^{-1} \\
& = gqg^{-1}(gqg^{-1}-(gqg^{-1})^*)^{-1}  . \end{array} 
$$
\end{num}

\begin{fed}\label{4.3} \rm 
We denote by $\prof =\{[a] \in \proj : pap \in G_{pAp}\} $.
The elements of $\prof$ are called the  ``finite points'' of $\proj $. In proposition
\ref{4.7} below we shall see that they are the analogue in this setting of the
finite points in the usual projective spaces and in the matrix projective spaces of 
Schwartz and Zaks.
\end{fed}

\begin{lem}\label{4.4}
 $\prof  \inc \proo = \{[up] : u \in \cu _A \}$.
Moreover, if $v \in \kp $ and $[v] \in \prof $, then $v \in \ua \cdot  p$.
\end{lem}

\noi Proof.  Let $v \in \kp $
such that $pvp = pv \in G_{pAp}$, which means that $[v] \in \prof$. 
We must prove  that  $v \in \cu _A\ p$. 
In matrix form, in terms of $p$, we can write $v = 
\bm v_1 & 0 \\ v_2& 0 \em $. Then 
$v_1 \in G_{pAp}$. Let $x = v_2 v_1^{-1} \in (1-p)Ap$. Then 
$$
v \ \sim _1 \bm p & 0 \\ x&0 \em  = p+x $$
The curve $\gamma (t) = [p+tx] $ joins
$[p]$ with $[v]$. Then $\varrho_p (v)=vv^* \in (\ep )_p$,
the connected component of $p$ in $\ep $. Since 
$(\ep )_p \inc \cu(p)$ by (\ref{2.16}), 
there exists $u \in \cu _A$ such
that $vv^* = upu^*$, i.e. $[v]=[up]$. Let $w \in
U(pAp)$ such that $v = upw$. Then
$$
v = upw = u\bm w & 0 \\ 0&1-p \em p
\in \cu _Ap.$$\QED

\begin{exa}\label{4.5}\rm 
Suppose that $A= B^{2 \times 2}$ for a \csta \
$B$ and $p = \bm 1&0\\ 0&0 \em $. 
We can embed $B$ into $\proj $ via
$$
B \ni b \mapsto \left[\bm 1&0\\ b&0 \em 
\right] \in \proj .
$$
Then
$$
[a]= \left[\bm a_{11}&0\\ a_{21}&0 \em 
\right] \in \prof  \quad \hbox{iff} \quad 
a_{11} \in G_B .
$$
It is also easy to see that the map $k : B \to \prof $ given by 
$$
k(b) = \left[ \bm 1&0\\ b&0 \em \right] \quad  , \quad b \in B ,
$$
is a homeomorphism (see (\ref{4.7}) below). 

\noi Let $g =\bm x&y\\ z&w \em \in \gla$. Denote by 
\begin{equation}\label{4.6}
D(g)=\{ b\in B \ : \ T_g(k(b)) \in \prof  \} \inc B.
\end{equation}
Consider the map $M_g : D(g)\to B$
defined by
$$
M_g(b) = k^{-1}(\ T_g \ (\ k(b)\ ) \ ) \quad ,\quad b \in D(g) 
$$
(the M\"oebius map on $B$ defined by
$g$). Easy calculations show that, for $b \in D(g)$, 
$$
M_g(b) = k^{-1} \left( \ \left[ \bm x+yb & 0 \\ 
z+wb&0 \em \right] \ \right) = (z+wb)(x+yb)^{-1},
$$
a picture that justifies the name M\"oebius map. Note
that 
$$
D(g) = \{ \ b \in B \ : \ x+yb \in G_B \ \}.
$$
This set is not easy to characterize, and could well be empty.
Another way to regard this domain is given in the following:
\end{exa}

\begin{pro}\label{4.7} Let $A$ be a \csta
\ and $p \in A$ a projection. Then
\ben
\item[i)] The linear manifold $H_p = (1-p)Ap$ is
C$^\infty$-diffeomorphic to $\prof $ via the map
$$
k : H_p \to \prof  \quad \hbox{ given by } \quad k(x) = 
\left[ \bm p & 0 \\ x&0 \em \right] = [p+x]
\quad , \quad x \in H_p.
$$
\item[ii)] The diffeomorphism $\varrho_p : \proj  \to \ep $
given by $\varrho_p([v]) = vv^*$ for $v \in \kp $
maps finite points onto projections $q$ such that
$\|p-q\|<1$. That is
$$
\begin{array}{rl}
\varrho_p(\prof ) &= \ \{ \ q \in U(p) \ : \ \|p-q\|<1 \ \} \\
   &= \ \{ \ q \in \ep  \ : \ d_r(p,q) < \pi /2 \ \}.
\end{array}
$$
\een
\end{pro}

\noi Proof.  ii) Let $v \in \kp $
such that $\|vv^* -p\|<1$. If $v = \bm v_1 & 0 \\ 
v_2&0 \em $, we have that $\|vv^*-p\|<1 \Rightarrow
\|v_1v_1^*-1\|_{pAp} <1 \Rightarrow v_1v_1^* \in G_{pAp}$.

On the other hand, $\|vv^* -p\|<1 $ implies that $ \|v_2^*v_2\|=
\|v_2v_2^*\|<1$. Because $v \in \kp $ we know that $v_2^*v_2+v_1^*v_1
= p$. Hence $\|v_1^*v_1-1\|_{pAp}=\|v_2^*v_2\|<1 $ and also
$v_1^*v_1 \in G_{pAp}$. Then $v_1 \in G_{pAp}$ and one
inclusion is proved.

Conversely, suppose that $v \in \kp $, $v=\bm v_1 & 0 \\ 
v_2&0 \em $ and $v_1 \in G_{pAp}$. Let
$x = v_2v_1^{-1}$ and $a =\bm p & 0 \\ 
x&0 \em = p+x \sim_1 v$. Let $a = w|a|$ be 
the polar decomposition of $a$. Note that $v \sim_1 a \sim_1 w$, and 
$$
|a|= (a^*a)^{1/2}= \bm (p+ x^*x)^{1/2} & 0 \\ 
0&0 \em \quad \Rightarrow \quad 
w = \bm (p+x^*x)^{-1/2} & 0 \\ 
z&0 \em .
$$
By  Lemma 4.4, we know that $w \in \cu _A \cdot p$.
Let $u \in \cu _A$ such that $w=up$. Then 
$$
u = \bm (p+x^*x)^{-1/2} & y \\ 
z&t \em .
$$
Now, 
$$
\|w^*w-p\|=\|upu^*-p\|=\|up-pu\| = 
\| \bm 0 & -y \\ 
z&0 \em \| .
$$
Since $u \in \cu _A$, $(p+x^*x)^{-1} +yy^* =
(p+x^*x)^{-1} + z^*z=p$.
Then
$$
\|yy^*\|=\|z^*z\| = \|p-(p+x^*x)^{-1}\|_{pAp}.
$$

\noi \bf Claim: \rm $ \|p-(p+x^*x)^{-1}\|_{pAp}<1$.

\noi Indeed, $p+x^*x \ge p \ \Rightarrow \ \sigma_{pAp}
((p+x^*x)^{-1}) \inc (0,1]$. Then $\sigma_{pAp}
((p+x^*x)^{-1}-1) \inc (-1,0] \ \Rightarrow \
\|(p+x^*x)^{-1}-1\|_{pAp}<1$, because it is selfadjoint.
Therefore $\|vv^*-p\|=\|ww^*-p\|= \max (\|y\| ,
\|z\|)<1 $, and the proof of (ii) is finished.
\smallskip
\noi
i) It is easy to see
that the map $k$ is bijective (note that
$[p+x]=[p+y]$ in $\proj  \ \Rightarrow \ x=y)$.
$k$ is C$^\infty$  because is the composition of the
C$^\infty$ maps $x \mapsto x+p$ and $x+p \mapsto [x+p]$.
Moreover, if $V_p = \{ q \in U(p) : \|p-q\|<1 \}$,
there exists a C$^\infty$ map (see \cite{[CPR6]})
$$
s_p : V_p \to \cu _A \quad \hbox{ such that } \quad 
s_p(q) \ p \ s_p(q)^* = q \quad , \quad q \in V_p.
$$
By ii) we know that $\varrho_p (\prof ) = V_p$. Then the map
$$\begin{array}{ccccccc}
     &  \varrho_p     &                &  s_p    &       && \\&\\
([p+x]) &  \mapsto  & q=\varrho_p ([p+x]) & \mapsto & s_p(q) &
\mapsto & (1-p) s_p(q)p \ [p \ s_p(q) \ p]^{-1} =x \\		  &\\
 && \zin       &&  \zin      && \zin      \\ 					 &\\
 &&  V_p    &&  \cu _A &&  H_p
\end{array}
$$
for $[p+x]\in \prof$, is the inverse of $k$ and  is
C$^\infty$. Note that 
$ps_p(q)p \in G_{pAp}$, since $[s_p(q)p] = [p+x] \in \prof $.\QED

\begin{rem}\label{4.8}
\rm Suppose that $A \inc {\cal B} (H)$ for a Hilbert 
space $H$. Via the identifications 
$$
k:H_p \to \prof \ \hbox{ and } \ 
\varrho_p : \prof \to \{q\in \ep : \|p-q\|<1\},
$$ 
of (\ref{4.7}), we can deduce that, for $g \in \gla$, the domain of 
the M\"oebius map $M_g$ induced by the projectivity $T_g$ is 
$$
D(g) = \{ q \in \ep : \|p-q\|<1 \ \hbox{ and } \ 
\|P_{g(p(H))}-p\|<1 \}.
$$
\end{rem}

\begin{num}\label{4.9}\rm 
In (\ref{4.1}) we defined the 
projectivities  $T_g$ for $g \in \gla$. Denote by 
$$
{\cal T} (\proj )\  = \ \{ \ T_g \ : \ g \in \gla \ \},
$$
the group of all projectivities on $\proj $. In order to characterize
the group ${\cal T }(\proj )$, we have to describe the isotropy group
$$
{\cal N} (\proj ) \ = \ \{ \ g \in \gla \ : \ T_g = Id_{\proj } \ \}.
$$
In \cite{[SZ2]} it is shown that, for $A=M_{2n}(\zC  )$,  ${\cal N} (\proj ) =
G_{\zC . I} $, the invertible scalar matrices. 
For a general \csta \ $A$, denote by 
$$
Z(A) \ = \ \{\ a \in A \ : \ ab = ba \ \hbox{ for all } \ b \in A \ \} 
$$
the center of $A$. 
\end{num}

\begin{pro}\label{4.10}
Let $A$ be a \csta \ and $p \in A$
a projection. Then
\ben
\item  ${\cal N}(\proj ) \ = \ \{ \ g \in \gla \ : \ g(q(H)) = q(H) \ 
\hbox{ for all } q \in \ep  \ \}$.
\item  $G_{Z(A)} \subseteq  \{ g \in \gla : 
gq=qg \quad
\hbox{ for all } \quad q \in \ep \ \} \subseteq {\cal N} (\proj )$.
\item  If $A = {\cal B} (H)$ for a separable Hilbert space $H$, then 
${\cal N}(\proj ) = G_{\zC  .I}$.
\item If $A$ is a von Neumann  factor of type III
(on a separable Hilbert space), 
then ${\cal N}(\proj ) = G_{\zC  .I}$.
\een
\end{pro}

\noi Proof. \rm Item 1 is apparent from the definitions. Since
$gq(H) = gqg^{-1} (H)$ for all $g \in \gla$ and $q \in \ep $, item 
2 follows from item 1. 

To prove item 3 consider first the case when 
$\dim p(H) = \infty$. In this case, for all $x \in H$ the 
subspace $<x>^\bot$ is the image of some
$q \in \ep ({\cal B} (H)) = \{$ projections $ q \in {\cal B} (H) \ : 
\ \dim q(H) = \infty \ \}$, 
since two projections are equivalent iff their images
have the same dimension in ${\cal B}(H)$. 
Therefore, if $g \in {\cal N}(\proj )$ then $x$ is an eigenvector for $g^*$
for all $x \in H$. Hence $g\in \zC  I$. 
If $\dim p(H) = n < \infty$ then any $g \in {\cal N}(\proj )$ should have all
subspaces of dimension $n$ as invariant spaces. It is easy to
see that such operators must be scalar multiples of the identity.

In order to prove item 4 recall that all 
non zero projections a factor of type III  are equivalent. If $ q \in A$ is a
projection such that $0 \ne q \ne 1$ then both $q \sim 1-q \sim p$. Let 
$g \in  {\cal N}(\proj )$. Then by item 1 we have that 
$$
qgq = gq \quad \hbox { and } \quad (1-q)g(1-q)=g(1-q).
$$
Hence $gq=qg$ for all projections $q \in A$, and $g \in Z(A) = \zC  I$.\QED

\begin{rem}\label{4.11}\rm 
We conjecture that for every \csta \ $A$ and $ p\in A$, 
$$
{\cal N} (\proj ) = \{ g \in \gla : gq=qg \quad
\hbox{ for all } \quad q \in \ep \ \} 
$$
Under reasonable hypothesis, this implies that ${\cal N} (\proj ) =
G_{Z(A)}$. A slight improvement of the argument used to show (\ref{4.10}), item 4
can be used to show that this conjecture is valid if $A$ is a type III von 
Neumann algebra  with separable predual and $p$ has 
central carrier $1$. 
\end{rem}

\section{The holomorphic structure of $\proj$ and $\ep$.}

As a homogeneous space of the complex analytic Lie group $G_A$, $\proj$
inherits a natural complex structure. It will be shown that the
projectivities $T_g$ are biholomorphic. These facts can be shown in an
explicit way, making use of the local charts $g\times\prof$, $g\in G_A$.

In \cite{[W]} Wilkins introduced a complex structure for the grassmannians. Under
the identification $\ep \simeq \proj$, both structures coincide.
Let us denote by 
$\varrho_q :\projq \to \eq = \ep$ the homeomorphisms (and isometries) 
of equation (\ref{2.13}), for all $q \in \ep$. 
\smallskip

\noi 
Let $A$ be a \csta \ and 
$p, q$ two equivalent projections in $A$. We consider the isometry
\begin{equation}\label{5.1}
\psi_{p,q}: \projq \to \proj \quad \hbox{ given by } \quad
\psi_{p,q} = \varrho_p^{-1} \circ \varrho_q .
\end{equation}

\begin{rem}\label{5.2}\rm 
By proposition \ref{4.7}, for each $q \in \ep$ we have that, 
\begin{equation}\label{5.3}
\psi_{p,q}(\profq) = \{ b \in \proj : d_c (b,\varrho_p^{-1}(q) )<1 \} :=
B_{\proj}(\varrho_p^{-1}(q),1), 
\end{equation}
because both sets are mapped onto $\{r \in \ep:\|r-q\|<1\}$ by 
$\varrho_q$ and $\varrho_p$, respectively. We denote 
by $ k_q:H_q = (1-q)Aq \to \profq $ the homeomorphisms 
of proposition \ref {4.7}, for each $q\in \ep$. 
We can define now the homeomorphisms
\begin{equation}\label{5.4}
k'_q : H_q \to B_{\proj} (\varrho_p^{-1}(q), 1) \quad \hbox{ given by }
\quad k'_q = \psi_{p,q} \circ k_q. \end{equation}
The maps $k'_q$, for $q \in \ep$ are almost the local charts for $\proj$. 
It just remains to uniformize the different Banach spaces $(1-q)Aq=H_q$, for 
different projections $q$. Note that the different connected components of
$\proj$ lie at chordal distance greater than $1$. Therefore in order to 
study the differential structure of $\proj$ we can work in each component. 
For simplicity, 
we shall define the complex structure only for the space $\proo $ which
is the union of several connected components of $\proj$. Note that 
$\varrho_p(\proo) = \{upu^* : u \in \cu_A\}=\cu (p)$, the unitary orbit of 
$p$. If  $q \in \up $ and  $w \in \cu_A$ such that $wpw^*=q$, then 
\begin{equation}\label{5.5}
Ad_w(H_p) = wH_pw^*=w(1-p)Apw^* = (1-q)Aq = H_q, 
\end{equation}
where $Ad_w$ is the inner automorphism of $A$ defined  by $w$.
\end{rem}

\medskip
\noi Let $a \in \proo$, $q=  \varrho_p(a)\in \up$ and $w\in \cu_A$ such that 
$wpw^*=q$. Using equations (\ref{5.3}), (\ref{5.4}) and (\ref{5.5}) we define the homeomorphism 

\begin{equation}\label{5.6}
\phi_a :H_p \to B_{\proj}(a,1) \quad 
\hbox{ given by } \quad \phi_a = k'_q \circ Ad_w 
=\psi_{p,q} \circ k_q \circ Ad_w.
\end{equation}

\begin{teo}\label{5.7} The family of local charts
$(\phi_a)_{a \in \proo}$ (choosing one appropiate $w$ for each a) defines
a complex holomorphic structure for the space $\proo$.
\end{teo}

\noi Proof. \rm We already know that all maps $\phi_a : H_p \to B_{\proj}(a,1)$
are homeomor\-phisms. So it remains to check taht these maps are  
compatibible with the analytic structure of $H_p$. 
In other words, if  $q, r \in \up$ and there exist $s \in \up $
such that $\|q-s\|<1$ and $\|r-s\|<1$, we must show that 
$(k'_q)^{-1} \circ k'_r$  is analytic. This will suffice because the maps
$Ad_w$ are analytic for all $w \in \cu_A$. 

\noi Case 1: Suppose that $\|q-r\|<1$. 
In this case, easy computations show the following formula: 
if $x \in H_r$ and 
$k'_r(x) \in  B_{\proj}(\varrho_p^{-1}(q), 1)= k'_q(H_q)$, then
\begin{equation}\label{5.8}
(k'_q)^{-1} \circ k'_r (x) = (1-q)(r+x)q. (q(r+x)q)^{-1}, 
\end{equation}
where the inverse of $q(r+x)q$ is taken in $qAq$. 
It is clear that the formula (\ref{5.8}) defines an analytic map of the 
variable $x$. 

\noi Case 2: Suppose that $q, r \in \up$ and there exist $s \in \up $
such that $\|q-s\|<1$ and $\|r-s\|<1$. Then, in the adequate domain, 
$$
(k'_q)^{-1} \circ k'_r = [(k'_q)^{-1}\circ k'_s]
\circ [(k'_s)^{-1}\circ k'_r].
$$
Since both maps on the right hand side  are analytic by Case 1,  the proof is
complete.\QED

\begin{rem}\label{5.`9}\rm The analytic structure can be 
extended to the whole  $\proj$ since, modulo the maps $\psi_{p,q}$, each connected
component of $\proj$ is included in 
$\prooq$ for some $q \in \ep$. Then the analytic  manifold structure can be 
defined around $q$ in the same way as in Theorem \ref{5.7}.
\end{rem}

\begin{rem}\label{5.10}\rm 
The following properties of $\proj$ are now easy to see:
\ben
\item  Each projectivity $T_g$, for $g \in \gla$, is biholomorphic.
\item The action of $\gla$ over $\proo$ given by the 
map $\pi_p : \gla \to \proo $ defined by $\pi_p(g)=T_g([p])$, 
$g \in \gla$, defines an analytic homogeneous space. The structure group is
the isotropy group 
$$
I_p = \{g \in \gla : T_g([p])=[p] \} = \{g \in \gla : 
(1-p)gp = 0 \ \hbox{ and } \ pgp \in G_{pAp} \},
$$ 
which is an union of connected components of the group of invertible 
elements of the subalgebra 
$$
T_p(A) = \{ a \in A : (1-p)ap =0 \} \inc A
$$ 
of $p$-upper triangular elements of $A$. This algebra is the tangent space 
at the identity of the group $I_p$. It is also the kernel of the 
differential $T(\pi_p)_1$ of $\pi_p$ at $1$, since 
$T(\pi_p)_1(a) = (1-p)ap$, for all $a \in A$.
\item The homogeneous space given by $\pi_p : \gla \to \proo $ admits a 
reductive structure given by the horizontal space $H_p =(1-p)Ap$ which can
be tranported homogeneously to all elements of $\gla$. Note that this 
horizontal space is precisely the domain of our local charts and can also be 
naturally identified with the tangent space 
$T(\proj)_{[p]}$ of $\proj$ at $[p]$.
\een
\end{rem}

\begin{rem}\label{5.11}\rm A complex structure can be 
defined in the Grassmannian $\up$ via  the map $\varrho_p$, i.e. pulling back 
the complex structure of $\proo$. This structure is compatible 
with the real structure, since $\varrho_p$ is a C$^\infty$ diffeomorphism 
by remark \ref{2.17}. It also allows 
us to define the analytic homogeneous reductive 
structure of $\up$ given by the new action of $\gla$ over $\up$:
$$
\pi_p : \gla \to \up \quad \hbox{ given by } \quad 
\pi_p(g) = P_{g(Im(p))} .
$$
Note that this action was described in (\ref{4.2}) as $\pi_p(g) = T_g(p)$. A
remarkable fact is that the formula for $T_g(p)$ given in (\ref{4.2}) becomes
analytic in the variable $p$, although the involution is 
involved in its description.
It can also be remarked that 
this complex structure agrees with the complex structure defined 
in the Grassmannians by Wilkins in \cite{[W]}. 
\end{rem}

\section{The non-Euclidean metrics on $\proj$.}

Suppose that $A$ is represented on a Hilbert space $H$. The projection $p$ induces
a Krein structure on $H$, by means of the selfadjoint symmetry 
$\eps=2p-1$. The set $\uep$ of operators of $A$ 
which are unitaries for this form
is called the group of $\eps$-unitaries. In this section we study
a subset $\d+ \subset \proj$, defined in (\ref{6.7,5}), 
which is homogeneous under the action of
$\uep$. Moreover, it is shown that $\d+$ can be regarded as a copy of
$\uep^+=\uep \cap A^+$ inside $\proj$, where $A^+$ denotes the space of
positive invertible elements of $A$. 
$\d+$ can be identified also with the space $N(A,p)$ of ``normal'' idempotents 
over the projection $p$, following the theory of \cite{[CPR5]} and \cite{[CPR6]}, pp 60. 

The space  $\uep^+$ is a totally geodesic
submanifold of $A^+$, which is a hyperbolic space, that is, a (non
riemannian) manifold of non positive curvature (in the sense of Gromov \cite{[G]}).
Therefore $\uep^+$ is a hyperbolic space in itself. In particular, it has
a rectifiable metric whose short curves can be explicitely computed. This
metric can be translated to $\d+$. This translation, which has a natural
intrinsic definition in terms of the already considered metrics of
$\proj$, will be called the  $non \ euclidean$   metric $E_n$. Summarizing,
$\d+$ will be shown to be a hyperbolic space inside $\proj$, with
an isometric action of $\uep$. 


\begin{fed}\label{6.2}\rm  \ 

\ben
\item Let $p \in A \inc {\cal B} (H)$ a projection. Consider the symmetry 
$$
\eps = 2p-1 = \bm 1&0\\0&-1 \em. 
$$
\item Denote by $\cu_\eps(A)$ the space of $\eps$-unitary elements of $A$, i.e. 
those $u \in A$ such that $<\eps u (\xi ), u (\eta )> 
= <\eps (\xi ),\eta>$, 
for all $\xi ,\eta  \in H$. Easy computations show that 
\begin{equation}\label{6.3}
\uep = \{ u \in A : u^*\eps u = \eps \} = 
\{u \in A : \eps u^*\eps = u^{-1}\}. 
\end{equation}
 \item Denote by 
\begin{equation}\label{6.4}
\uep^+  = \uep \cap A^+, \end{equation}
the set of positive $\eps$-unitary elements of $A$.
\een
\end{fed}

\smallskip
\noi In the following remark we state several well known properties 
of the sets $\uep$ and $\uep^+$ (see \cite{[CPR5]} and \cite{[CPR6]}) :

\begin{rem}\label{6.5} \rm \  

\ben
\item $\uep$ is a closed subgroup of $\gla$. Actually 
it is a real Banach-Lie group.
\item If $u \in \uep$, then $u^*, \ u^*u $ and $u^{-1} \in \uep$. 
\item  $\lambda \in \uep^+$ if and only if 
$\lambda \eps = \eps \lambda ^{-1}$. 
\item For all $\lambda \in \uep^+ $,  $ X= \log \lambda \in A $, 
verifies that $X=X^*$ and $\eps X = -X\eps$. 
\item In matrix form, we have that
$\lambda \in \uep^+ $ if and only if there exists   $x \in H_p$ 
such that 
$$
\lambda = e^X \quad \quad \hbox{ where } \quad X = {\bm 0&x^*\\x&0 \em }.
$$
In this case, $x$ is \bf unique  \rm.
\item Using item 5, one deduces that 
$$
\lambda \in \uep^+ \Rightarrow 
\lambda^t \in \uep^+ \quad \hbox{ for all } \quad t\in\zR. 
$$
\item In particular, if $u \in \uep$, then 
$|u|=(u^*u)^{1/2} \in \uep$. In other words, the unitary and positive parts of
each $u\in \uep$ in its polar decomposition remain in $\uep$. Note also that 
$\ua \cap \uep = \{u \in \ua : u\eps = \eps u \} =\{u \in \ua : up = p u \}$.

\item The metric of $A^+$ and its geodesics (see \cite{[CPR4]}, \cite{[CPR5]} and 
\cite{[CPR7]}) can be restricted to $\uep^+$.
Indeed, if $\lambda ,  \ \mu \in A^+$, the unique 
geodesic of $A^+$ joining them is given by 
$$
\gamma_{\lambda  \mu } (t) = 
\mu ^{1/2}(\mu ^{-1/2}\lambda\mu ^{-1/2})^t\mu ^{1/2}
\ , \quad t \in [0,1] .
$$
Using items 1 and 6 one shows that if $\lambda ,  \ \mu \in \uep^+$ then
$\gamma_{\lambda  \mu } (t)\in \uep^+$. 
\item The Finsler structure of $A^+$ (see \cite{[CPR5]}) induces  
a rectifiable metric on $A^+$ given by
$$
d_+(\lambda , \mu ) = \|\log (\mu ^{-1/2}\lambda\mu ^{-1/2}) \| \ \ 
\hbox{ for } \ \ \lambda, \ \mu \in A^+.
$$
For this metric the geodesics $\gamma_{\lambda  \mu }$ are of minimal length. 
Restricted to $\uep^+$, this metric is also rectifiable by item 8, because the
geodesic curves remain in $\uep^+$ and are of course of minimal length. 
\een
\end{rem}

\begin{rem}\label{6.6}\rm 
Using item 5 of (\ref{6.5}), easy computations, 
very similar to those of Theorem 3.5, show that for each 
$\lambda \in \uep^+$ there exists a unique $x \in H_p$ such that
$$
\lambda =\bm \cosh(|x|) &x^*
(\frac{\sinh t}{t}) (|x^*|)\\
x (\frac{\sinh t}{t}) (|x|) & \cosh (|x^*|) \em .
$$
In particular, using item 7 of (\ref{6.5}), this implies that for all $u \in \uep$, 
$$
\|(1-p)up\|=\|(1-p)u^{-1}p\| = \|(1-p)u^*p\|= \|\sin (|x|)\|,
$$
where $x\in H_p$ verifies that $|u| = e^{x+x^*}$ 
\end{rem}

\begin{num}\label{6.7}\rm 
Return now to the projective space $\proj$. Given $u \in A$, easy 
matrix computations using (\ref{6.3}) show that if $u = 
\bm a & c 
\\b&d\em $ then $u \in \uep$ if and only if $u$ is 
invertible and the following three conditions hold:
\ben
\item[i)] $a^*a-b^*b=p$,
\item[ii)] $d^*d -c^*c = 1-p$ and 
\item[iii)] $a^*c - b^*d = 0$.
\een
\end{num}

\medskip 
\noi Consider the set
\begin{equation}\label{6.7,5}
\d+ = \{[up] \     : \ u \in \uep \} \inc \proj .
\end{equation}
Denote by 
\begin{equation}\label{6.8}
K'_p(A) = \{ \bm a &  0
\\b&0\em \in Ap : a \in G_{pAp} \ \hbox{ and } \ 
a^*a-b^*b = p  \} \inc \lp .
\end{equation}

\begin{pro}\label{6.9} 
$$
\begin{array}{rl}
\d+ & = \{[u] : u \in K'_p(A)\} \\&\\
&     = \{[\lambda p ] : \lambda \in \uep^+ \}\\&\\
& =  \{\left[ \bm a &  0
\\b&0\em \right] : a \in G_{pAp} \ \hbox{ and } \ 
a^*a-b^*b > 0  \ \hbox{ in } \ 
pAp \}      \inc \prof , \end{array}
$$
and the map $Y : \uep^+ \to \d+ $ given by $Y(\lambda) = 
[\lambda^{1/2} p]$ is a
homeomorphism.
\end{pro}

\noi Proof. \rm We have the following trivial inclusions:
$$\begin{array}{rl}
 \{[\lambda p ] : \lambda \in \uep^+ \} & \inc   \d+ \\&\\
& \inc \{[u] : u \in K'_p(A)\} \\&\\
& \inc  \Big{ \{ } \left[\bm a &  0
\\b&0\em \right] : a \in G_{pAp} \ \hbox{ and } 
a^*a-b^*b > 0  \ \hbox{ in } \ 
pAp \Big{\}} . \end{array}
$$
So we have to show that 
$ \{ \left[\bm a &  0
\\b&0\em \right] : a \in G_{pAp} \ \hbox{ and } \ 
a^*a-b^*b > 0  \ \hbox{ in } \ 
pAp \}\inc  \{[\lambda p ] : \lambda \in \uep^+ \}$. Let 
 $v=  \bm x &  0
\\y&0\em $ such that $x \in G_{pAp}$ and 
$d = x^*x -y^*y >0$. Then 
$w=vd^{-1/2} \in K'_p(A)$, $[v]= [w]$  and $w = \bm a &  0
\\b&0\em $ with $a^*a-b^*b = p$. 
Since $a \in G_{pAp}$, by taking 
the (right) polar decomposition of $a$ in $pAp$ we 
can also suppose that $a>0$, that is $a = (p+b^*b)^{1/2}$.
Consider
\begin{equation}\label{6.10}
\lambda = \bm a &  b^* \\
b&  (1-p + bb^*)^{1/2}    \em > 
\bm |b| &  b^* \\ b&  |b^*|   \em \ge 0 , 
\end{equation}
because $a = (p+b^*b)^{1/2} > (b^*b)^{1/2} = |b|$ in $pAp$ and 
$(1-p+bb^*)^{1/2} > |b^*|$ in $(1-p)A(1-p)$. 
Then $\lambda \in A^+$. It is easy to see that $\lambda$ verifies 
the three conditions of (\ref{6.7}). So $\lambda \in \uep^+$ and $[\lambda p]=[v]$.

The map $Y_0:\uep^+ \to \d+$ given by $Y_0 (\lambda ) = [\lambda p]$ is 
therefore continuous and surjective. To see that $Y_0$ is
injective, suppose $\lambda , \mu \in \uep^+$ with $[\lambda p ] = [\mu p]$. 
Put   $\lambda p = a+b$ and $\mu p = c + d $ with $a,c \in pAp^+$ and 
$ b,d \in H_p$. Then 
$$
ba^{-1} = b(p+b^*b)^{-1/2} = d(p+d^*d)^{-1/2} = dc^{-1}.
$$
Taking their polar decompositions in ${\cal B} (H)$, both elements 
have the same 
partial isometry,  say $u$, and therefore 
$$
ba^{-1} = u|b|(p+|b|^2)^{-1/2} = u|d|(p+|d|^2)^{-1/2} = dc^{-1},
$$
proving that $u$ is also the partial isometry for $b$ and $d$ in 
their polar decompositions. This 
implies that $|b| = |d|$ since that map $f(t) = \frac{t}{(1+t^2)^{1/2}}$
has inverse $g(s)=\frac{s}{(1-s^2)^{1/2}}$. Then $b=d$ and $\lambda = \mu$ 
by (\ref{6.7}) and equation (\ref{6.10}). Note that we have already constructed the inverse of 
$Y_0$ by passing through $H_p$: 
$$ 
Y_0^{-1}(\left[ \bm a &  0 \\b&0\em \right] ) = 
\bm (p+d^*d)^{1/2} &  d^*
\\d&(1-p+dd^*)^{1/2} \em \in \uep^+ ,
$$
where $d = ba^{-1} (1-|ba^{-1}|^2)^{-1/2}$. Clearly this map
is also continuous. Finally, note that the map $Y$ is the composition of the 
homeomorphism of $\uep^+$ which consists of taking square roots, 
with the homeomorphism $Y_0$. 
Then the proof is complete.\QED

\begin{rem}\label{6.11} \rm \ 
\ben
\item From the proof of (\ref{6.9}), it follows quite easily that
$$
K'_p(A) = \uep \cdot p .
$$
\item
We have shown some characterizations of $\d+$ in terms of its representatives 
in $A$. Now we give other characterizations of $\d+$ in terms of the three
natural metrics on $\prof$: the chordal and spherical metrics and the new
metric $d_k$ on $\prof$ given by the map $k_p^{-1}:\prof \to H_p $ of (\ref{4.7}): 
$$
d_k(l,m) = \|k_p^{-1}(l)-k_p^{-1}(m)\| \quad \hbox{ for } \quad 
l,\ m\in \prof . 
$$
\een
\end{rem}

\begin{cor}\label{6.12} 
$$
\begin{array}{rl}
\d+ & = \{m \in \prof : d_k(m, [p]) = \|k_p^{-1}(m)\|  <1 \} \\&\\
& = \{m \in \prof : d_c(m,[p]) = \|\varrho_p(m) -p\| < 
\frac{\sqrt{2}}{2} \} \\&\\
    & = \{m \in \prof : d_r(m , [p]) < \pi /4 \}. \end{array}
$$
\end{cor}

\noi Proof. \rm Let $m\in \d+$ and $\lambda \in \uep^+$ such that 
$[\lambda p ] = m $. If 
$\lambda p =  \bm a &  0 \\b&0 \em
$ then 
$$
\|k_p^{-1}(m)- k_p^{-1}([p])\| = \| ba^{-1} \| = 
\|b(p+b^*b)^{-1/2}\| < 1.
$$ 
On the other hand, if $x \in H_p$ and $\|x\|<1$ then
$k_p(x) = [p+x] $ and $p-x^*x >0 $ in $pAp$. Then $k_p(x) \in \d+$.

\noi 
In order to prove the other equalities, recall from the proof of (\ref{3.5}), that 
if $m \in \prof$, then there exist $y \in H_p$ such that $d_r(m,[p])=\|y\|$, 
and 
$$
m = \left[\bm \cos(|y|) &  0 \\
y (\frac{\sin t}{t}) (|y|) & 0 \em \right]
\quad \Rightarrow \quad k_p^{-1}(m) = y (\frac{\tan t}{t}) (|y|) .
$$
Hence  
\begin{equation}\label{6.13}
d_k(m,[p]) = \|k_p^{-1}(m)\|= 
\|y(\frac{\tan t}{t})(|y|)\|=\|\tan (|y|)\| = \tan (\|y\|), 
\end{equation}
because the map $f(t) = \tan (t) $ 
is monotone increasing. Therefore $\|k_p^{-1}(m)\|<1$ if
and only if $d_r(m,[p]) = \|y\|= \arctan (d_k(m,[p])<\pi /4$. 
Now the remainding equality becomes apparent using that, by (\ref{3.5}), 
$d_c(m,[p]) = \sin (d_r(m,[p]))$.\QED

\begin{rem}\label{6.14} \rm From (\ref{3.5}) and the proof of (\ref{6.12}), 
we can deduce the following facts: for all $m,n \in \prof$,
\ben
\item  $d_k(m,[p]) = \tan (d_r(m,[p]))$.
\item  $d_c(m,n) = \sin (d_r(m,n))$.
\een
\end{rem}
\medskip

\noi
Note that the three metrics have clear geometrical senses: 
the chordal metric
is the one associated to the norm in $\cu (p)$ via the identification 
(which is in fact a bi-analytic map) $\varrho_p$. The metric $d_r$  is 
the rectifiable metric generated by $d_c$ taking the infima of 
the lengths of curves (and having the geodesics of the linear connection 
as minimal curves). On the 
other hand $d_k$ is the metric induced on $\prof$ by the atlas of 
local charts of its complex manifold structure. Note also that 
they are related
by the previous formulae and  depend on the norm of some particular
vectors in $H_p$, which is the tangent space at $p$ of $\prof$.
On $\d+$ we have a fourth metric, induced by the metric $d_+$ (see (7) of 
(\ref{6.5})) of $\uep^+$ via the map $Y: \uep^+ \to \prof$ of (\ref{6.9}) 
(which is a C$^\infty$ diffeomorphism). Let $m,n \in \d+$ and 
$\mu , \nu \in \uep^+$ such that $m=[\mu^{1/2} p]=Y(\mu)$ and 
$n = [\nu^{1/2} p ]=Y(\nu)$. Then
$$
d_+(m,n) = d_+(\mu , \nu) = \|\log (\nu^{-1/2} \mu \nu^{-1/2})\|.
$$ 
\medskip
\noi
We define now the non-Euclidean metrics on $\d+$ following \cite{[SZ2]}:

\begin{num}\label{6.15}
\rm Let $m,n \in \d+$ and 
$u,v \in \uep $ such that $m=[up]$ and $n=[vp]$. We consider the following
three functions:
\ben
\item $\rho (m,n) = \|(1-p) u^* \eps v p\|$. This function is clearly well 
defined but is not a metric. 
\item The ``pseudo-chordal'' metric:
$$ d_{pc} (m,n) = \frac{\rho (m,n)}{(1+\rho (m,n)^2)^{1/2}}. $$
\item The non-Euyclidean metric:
$$
E_n(m,n)= \frac12 \ \log \ \frac {1+d_{pc} (m,n)}{1-d_{pc} (m,n)}.
$$
\een
\end{num}

\begin{rem}\label{6.16}\rm 
In order to relate the metrics just defined, we recall from 
item 5 of (\ref{6.5}) the action of the 
group $\uep$ over $\uep^+$ given by $u \times \lambda = u \lambda u^*$, for 
$u \in \uep $ and $\lambda \in \uep^+$. By item 8 of (\ref{6.5}), this
action is isometric with respect to the Finsler metric of $\uep^+$, since the same action
is isometric at $A^+$ (see \cite{[CPR5]}). 
Therefore this action is also isometric for the geodesic 
metric $d_+$ defined in item 9 of (\ref{6.5}). This fact yields
$$
d_+(\mu , \nu ) = d_+(\nu ^{-1/2} \mu \nu^{-1/2} , 1) = 
\|\log (\nu ^{-1/2} \mu \nu^{-1/2})\|. 
$$
We consider also the action of $\uep$ over $\d+$ induced via the map $Y$ of 
(\ref{6.9}). More explicitely,  for $u \in \uep $ and $\lambda \in \uep^+$, 
$$
u\times [\lambda^{1/2} p] = u \times Y(\lambda) := Y(u \times \lambda) =
[(u\lambda u^*)^{1/2} p].
$$ 
\end{rem}

\begin{pro}\label{6.17} \

\ben
\item The pseudo-chordal and non-Euclidean metrics are symmetric and 
invariant under  the action of $\uep$ on $\d+$  defined in (\ref{6.16}).  
\item  For  $m,n \in \d+$, let $Y^{-1}(m) = \mu \in \uep^+$ and
$Y^{-1}(n) = \nu \in \uep^+$. Then   
$$
d_{pc}(m,n) =  \|k_p^{-1}([\nu ^{-1/2} \mu ^{1/2} p ])\|=
\|k_p^{-1}([ \nu ^{1/2} \mu^{-1/2} p])\| \ \ \hbox{ and, in particular } 
$$
$$
d_{pc}(m, [p]) = d_k (m, [p]). 
$$ 
\een
\end{pro}

\noi Proof. \rm 1) First note that, using (\ref{6.6}), 
$$
\rho (m,n) =\|(1-p) \mu \eps \nu p \| =
\|(1-p)\mu \nu^{-1} p\| = \|(1-p)\nu \mu^{-1} p \| = \rho(n,m), 
$$
and $\rho$ is symmetric. On the other hand, if  $ u \in \uep $, its left 
polar decomposition is given by $u = (uu^*)^{1/2} w$, where 
$w\in \ua \cap \uep$ commutes  with $p$ by item 7 of (\ref{6.5}). Hence, for all 
$u  \in \uep$ we have that 
\begin{equation}\label{6.188}
[up] = [(uu^*)^{1/2} w p] = [(uu^*)^{1/2} p (pwp)] = [(uu^*)^{1/2}p]
=[|u^*|p].
\end{equation}
Now we can describe more clearly the action of $\uep $ over $\d+$ of (\ref{6.16}): 
let $u \in \uep $ and $\lambda \in \uep^+$, then, by (\ref{6.18}), 
$$
u \times [\lambda ^{1/2}p] = [(u\lambda u^*)^{1/2} p] =
[u\lambda^{1/2} p].
$$
Then, for $u \in \uep$, 
$$
\begin{array}{rl}
\rho(u\times [ \mu^{1/2} p],u\times [ \nu^{1/2}p]) 
& = \rho([u\mu^{1/2} p], [u\nu^{1/2} p]) \\&\\
& = \|(1-p) \mu^{1/2} u^* \eps u \nu^{1/2}  p \| \\&\\
& = \|(1-p) \mu^{1/2} \eps  \nu^{1/2}  p \| \\&\\
& = \rho( [ \mu^{1/2} p],[\nu^{1/2} p]) . \end{array}
$$
Therefore $\rho$ is symmetric and invariant under the action of 
$\uep$ on $\d+$. It is clear that the same happens for $d_{pc}$ 
and $E_n$, since they are defined in terms of $\rho$.

\noi 2) Using 1)  we have that 
$$
\begin{array}{rl}
d_{pc}(m,n) & = d_{pc}([\mu^{1/2}p],[\nu^{1/2}p]) \\&\\
            & = d_{pc}([\nu^{-1/2}\mu^{1/2}p],[p]) \\&\\
            & = d_{pc}([(\nu^{-1/2}\mu \nu^{-1/2})^{1/2}p], [p]). 
\end{array}
$$
Then, using again (\ref{6.16}), we can suppose that $n = [p]$, since the action 
of $\nu^{-1/2}$ transforms $n= [\nu^{1/2}p]$ to $[p]$ in $\d+$ and 
$\nu$ to $1$ in $\uep^+$. Now, let $\mu^{1/2} p = 
 \bm a &  0 \\b&0 \em $. Then $x= k_p^{-1}(m) = 
ba^{-1} = b(p+b^*b)^{-1/2}$ and 
$$
d_k(m, [p])^2 = \|x^*x\| = \|b^*b(p+b^*b)^{-1}\| = \frac{\|b\|^2}{(1+\|b\|^2)}.
$$
Note that $\|b\|= \|(1-p)\mu^{1/2} p\| = \rho([\mu^{1/2} p],[p]) =
\rho (m,[p])$. Therefore
$d_{pc}(m,[p]) = d_k(m,[p])$. \QED

\begin{teo}\label{6.18} \ 

\ben
\item If  $m,n \in \d+$, let $Y^{-1}(m) = \mu \in \uep^+$ and
$Y^{-1}(n) = \nu \in \uep^+$. Then   
$$
2 \ E_n(m,n) = d_+(m,n) = d_+(\mu ,\nu ) = \|\log (\nu^{-1/2}\mu\nu^{-1/2})\|,
$$
where $d_k$ is the metric on $\prof $ defined in (\ref{6.11}) and $d_+$ is the 
metric on $\d+$  defined in  (\ref{6.14}). 

\item The map $Y:\uep^+ \to \d+$ 
allows the translation of  
the C$^\infty$ homogeneous reductive  structure and Finsler metric 
of $\uep^+$ to $\d+$. In this sense, the 
geodesics $ \gamma _{\mu,\nu}$ defined in item 8 of (\ref{6.5}) 
yield minimal length geodesics in $\d+$ via $Y$: 
$$
 E_n(m,n) = \frac 12 \ d_+(\mu ,\nu) = 
\frac 12 \ \ell_{\uep^+} (\gamma _{\mu,\nu}) =
\  \ell_{\d+} (\gamma _{m,n}),
$$
where 
$$
\begin{array}{rl}
\gamma _{m,n} (t) & = Y\circ \gamma _{\mu,\nu}(t) \\ &\\
   & = Y(\mu ^{1/2}(\mu ^{-1/2}\nu\mu ^{-1/2})^t\mu ^{1/2}) \\&\\
 & = [\mu^{1/2}(\mu^{-1/2}\nu^{1/2}\mu^{-1/2})^{t/2} p]\ , 
 \quad t \in \zR .
\end{array}
$$
\een
\end{teo}

\noi Proof. \rm
Using (\ref{6.6}), there exists $z \in H_p$ such that, 
if $Z =    \bm 0 &  z \\z^*&0\em$, then
$$
\mu^{1/2} = e^{Z } =
\bm \cosh(|z|) &z^*(\frac{\sinh t}{t}) (|z^*|)\\
z (\frac{\sinh t}{t}) (|z|) & \cosh (|z^*|) \em
$$
and, by (\ref{6.5}),  
$$
d_+(m, [p]) = d_+ (\mu, 1)=\|\log(\mu)\|=  2   \|\log (\mu^{1/2}) \|= 
2 \| \bm 0 &  z \\z^*&0\em \| = 2 \|z\|.
$$
Easy computations similar as those of (\ref{6.13}), show that
$$ 
d_{pc}(m,[p]) = d_k(m,[p]) = \|z(\frac{\tanh (t)}{t}) (|z|)\|=
\|\tanh(|z|)\| = \tanh (\|z\|).
$$
Therefore $d_+(m,[p]) = 2$ argtanh$(d_{pc}(m,[p]))$. 
Elementary computations 
show that, for all $t\in (-1,1)$, argtanh$ (t) = \frac12 \ \log \ 
\frac{1+t^2}{1-t^2}$. Therefore we have proved that the metrics 
$2 E_n$ and $d_+$ coincide.

\noi 2) It follows because the geodesics $\gamma_{\mu, \nu}$ are 
minimal for $d_+$ in $\uep^+$ and we have translated the metric and Finsler 
structure from $\uep^+ $ to $\d+$ via $Y$. \QED


\vglue1truecm

\noindent{Esteban Andruchow}

\noindent{Instituto de Ciencias, UNGS,  San Miguel, Argentina}


\noindent{e-mail : eandruch@ungs.edu.ar}
\bigskip

\noindent{Gustavo Corach}

\noindent{Depto. de Matem\'atica, FCEN-UBA, Buenos Aires, Argentina and }

\noindent{Instituto Argentino de Matem\'atica, Buenos Aires,  Argentina}

\noindent{e-mail: gcorach@mate.dm.uba.ar}

\noindent{and}

\noindent{Demetrio Stojanoff}

\noindent{Depto. de Matem\'atica, FCE-UNLP,  La Plata,  Argentina and }

\noindent{Instituto Argentino de Matem\'atica, Buenos Aires,  Argentina}

\noindent{e-mail: demetrio@mate.dm.uba.ar}
\bigskip
 
\noindent{Mathematical Subject Classification 1991: Primary 46L05, 58B20.}

\end{document}